\documentclass{article}    

\usepackage{amsmath,amssymb,amsthm,epsfig,url,xcolor}      

\newtheorem{theorem}{Theorem}[section]
\newtheorem*{theorem*}{Theorem}
\newtheorem{proposition}[theorem]{Proposition}
\newtheorem{lemma}[theorem]{Lemma}

\newtheorem{corollary}[theorem]{Corollary}
\newtheorem*{corollary*}{Corollary}

\newtheorem*{conjecture*}{Conjecture}

\theoremstyle{remark}

\theoremstyle{definition}

\newcommand{\hide}[1]{{}}  

\newcommand{\ES}{Erd\H{o}s-Selfridge }
\newcommand{\vb}{{\bf v}}
\newcommand{\bfour}{\begin{tabular}{c|c|c|c}}
\newcommand{\efour}{\end{tabular}}
\newcommand{\bthree}{\begin{tabular}{c|c|c}}
\newcommand{\ethree}{\end{tabular}}
\newcommand{\bfive}{\begin{tabular}{c|c|c|c|c}}
\newcommand{\efive}{\end{tabular}}
\newcommand{\bsix}{\begin{tabular}{c|c|c|c|c|c}}
\newcommand{\esix}{\end{tabular}}
\newcommand{\bseven}{\begin{tabular}{c|c|c|c|c|c|c}}
\newcommand{\eseven}{\end{tabular}}

\newcommand{\HH}{\mathcal{H}}
\newcommand{\Hempty}{\mathcal{H}_{\mathrm{empty}}}

\begin{document}

\title{Cubic tic-tac-toe:  A matching-based approach}

\author{
John W.~Cain\footnote{Department 
of Mathematics, Harvard University, One Oxford Street, Cambridge, Massachusetts 02138, USA.
\newline
Corresponding author's e-mail: jcain2@math.harvard.edu},
Ioannis M.~Raymond, and Nora C.~K\"{a}llersj\"{o}
}

\maketitle

\bigskip

{\small
\begin{center}
{\sc Abstract}
\end{center}
In the natural generalization of tic-tac-toe to an $n \times n \times n$
board where $n \in \mathbb{N}$, it is known that the first player has a winning strategy if $n \leq 4$
and that either player can force a draw if $n \geq 8$.  The question of whether
the first player has a winning strategy if $n  = 5, 6$ or $7$ has remained open.
Here, we prove that the first player does not have a winning strategy if $n = 7$.
The proof, which is computer-assisted, exploits the fact that the second player's first four moves can 
always be chosen in such a way that their remaining moves can be automated via a 
simple pairing strategy.  The process of finding the pairing strategy involves
reframing the problem in such a way that the goal is to seek a maximal matching in
a bipartite graph that represents the tic-tac-toe board after each player has made
four moves.  We use the Hopcroft-Karp matching algorithm to find such maximal matchings.

\vspace{0.5 true cm}
}

\bigskip

\large


\section{Background} 

Nearly half a century has passed since Oren Patashnik published a proof that the first player has a winning 
strategy in $4 \times 4 \times 4$ tic-tac-toe~\cite{patashnik}. Patashnik's result represented one of the first 
major victories within the realm of {\em computer-assisted proof}. Just three years earlier, Appel and 
Haken~\cite{appel-H} had presented a proof of the four-color theorem, one of the most famous computer-assisted 
proofs of all time.  Those were exciting times---research computing clusters were becoming considerably more 
powerful, and personal computers were beginning to hit the market.

Patashnik's work, published in 1980, is actually the most recent significant advancement in our understanding
of strategies and outcomes in cubic tic-tac-toe games.  Since then, it has been known that the first player
has a winning strategy in $n \times n \times n$ tic-tac-toe if $n \leq 4$, and that either player can force
a draw (tie) if $n \geq 8$.  The cases $n = 5, 6, 7$ have remained unresolved, but here we prove that either  
player can force a draw if $n = 7$.  Our proof is computer-assisted, but very different than the
pruning-the-game-tree approach described in~\cite{patashnik}.  Rather, we will explain how a basic idea from
introductory graph theory (finding matchings in bipartite graphs) can be leveraged to prove that the second
player in $7 \times 7 \times 7$ tic-tac-toe can prevent the first player from winning by carefully pairing
each move with the immediately preceding move of their opponent.  Patashnik remarked that his exploration of
the $4 \times 4 \times 4$ game required 1500 hours of computer time, which would have cost about 50 million
dollars at then-commercial rates if the Yale University Department of Computer Science had not generously allowed
access to their computing cluster.  By contrast, our computational casework for the $7 \times 7 \times 7$  
game required only about five hours on a standard laptop computer.

In spite of the seemingly playful nature of this article's subject matter, there are several reasons why
we believe that the methods described herein are important.  First, generalized tic-tac-toe has
connections to major theorems in combinatorics.  For example, given any fixed $n \in \mathbb{N}$,
the Hales-Jewett theorem explains why the $n^d$ game with standard rules (described below) cannot end in a
draw if $d$ is sufficiently large.  In our estimation, devoting effort towards some of the {\em many}
open problems related to tic-tac-toe seems likely to generate results of independent importance in
graph theory and combinatorics.  Second, the proof of our main result offers a worked example of how one might
frame computational problems using the language of {\em matchings} and matching algorithms.  There are
many practical problems for which it is desirable to find maximal matchings efficiently, often involving
scheduling or assignment (e.g., matching medical residency applicants with hospitals based upon their
preference lists).

In the interest of keeping this manuscript reasonably self-contained, much of what follows is a survey article. We 
begin with an overview of tic-tac-toe, its variants and generalizations to higher dimensions, and some known 
results (Section~\ref{section2}).  An excellent and enjoyable-to-read reference for this material is the text of 
Beck~\cite{beck}, which also includes a compendium of related open problems.  Section~\ref{section44} 
establishes a framework for our main result, using an easily understood test case:  tic-tac-toe on a
$4 \times 4$ game board.  In doing so, we provide a worked example of (i) how to represent the state of the
game board using hypergraphs; (ii) how to describe defensive pairing strategies using the language of
matching theory; and (iii) how to actually find matchings that a player can leverage in order to prevent their
opponent from winning.  Readers familiar with combinatorial games are invited to skip directly to 
Section~\ref{section-main}, the proof of our main result concerning the $7 \times 7 \times 7$ game.  
Our proof is constructive, in that it explains how either player can prevent the other from winning.
Section~\ref{section-bigamy} outlines our attempt to provide a more elegant, non-constructive proof by
invoking a variant of the K\H{o}nig, Hall, and Egerv\'{a}ry matching theorem.  Finally, in 
Section~\ref{section-discussion} we discuss the two remaining unresolved cubic
tic-tac-toe games:  $5 \times 5 \times 5$ and $6 \times 6 \times 6$.  We also survey a very 
different approach for determining whether a player has a winning strategy, namely
the use of danger potentials as opposed to matching theory.


\section{Tic-Tac-Toe:  Rules, Game Play, and Known Results}  \label{section2}
  
Without further ado, let us explain the game of tic-tac-toe, its standard rules, and one useful variant of the
standard rules.

\subsection{What is Tic-Tac-Toe?}  \phantom{hello}

{\em Tic-tac-toe} is a two-player positional game known by a variety of names around the world (e.g., {\em
noughts-and-crosses}).  The two players take turns occupying positions on the game board, using the letters
X (first player) and O (second player) to indicate which positions they have occupied.  Each player has
some winning objective that they seek to achieve and, whichever player (if any) is first to achieve their
winning objective is declared the winner.  If neither player achieves their winning objective before all
positions in the game board have been occupied, the game is a {\em draw} or {\em tie}.  The most commonly
used board is a $3 \times 3$ grid (three rows and three columns) containing nine positions as sketched below.
$$
    \bthree
      &   &   \\ \hline
      &   &   \\ \hline
      &   &
  \ethree
$$
Under ``standard" rules (described below) on a $3 \times 3$ game board, both players have the same winning
objective:  occupy all three positions along one of eight {\em winning lines} (a row, a column, or one
of the two ``main" diagonals).


\subsection{Game Board and Winning Lines}  
  
Here is a more careful description of the tic-tac-toe game board and winning lines; rules will be discussed
afterwards.  Let $n$ and $d$ be positive integers, and let $[n]$ and $[d]$ denote the sets
$\{1, 2, \dots, n\}$ and $\{1, 2, \dots, d\}$, respectively.  The game of {\em $n^d$ tic-tac-toe} uses the
following {\em game board}:
$$
 B = \underbrace{[n] \times [n] \times \cdots \times [n]}_{d \mbox{ times}}.
$$
\noindent
Each {\em position} or {\em cell} in the board can be represented by a $d$-tuple of the form
$\vb = (a_{1}, a_{2}, \cdots, a_{d}) \in B$.  By a {\em winning line}, we mean an ordered list of $n$
such $d$-tuples $( \vb^{(1)}, \vb^{(2)}, \cdots, \vb^{(n)} )$ with the following two properties:
\begin{itemize}
  \item For each $i \in [d]$, the sequence $\vb_{i}^{(1)}, \vb_{i}^{(2)}, \dots, \vb_{i}^{(n)}$
        either increases from $1$ to $n$, decreases from $n$ to $1$, or remains constant.  Here,
        $\vb_{i}^{(j)}$ denotes the $i$th component of $\vb^{(j)}$.
  \item For at least one $i \in [d]$, the sequence $\vb_{i}^{(1)}, \vb_{i}^{(2)}, \dots, \vb_{i}^{(n)}$
        is non-constant.  
\end{itemize}
For example, in the $3^2$ game,
$$
 ( (2,1), (2,2), (2,3) )
$$
is an example of a winning line.  If one interprets $(i, j)$ as representing the position in row $i$ 
and column $j$ in a $3 \times 3$ game board, this winning line corresponds to the middle column.
In the $7^3$ game,
$$
(  (7,5,1), (6,5,2), (5,5,3), (4,5,4), (3,5,5), (2,5,6), (1,5,7) )
$$
is an example of a winning line.  The first coordinate decreases from 7 to 1, the second coordinate
is fixed, and the third coordinate increases from 1 to 7.  Enumerating the winning lines is crucial, so
we will make heavy use of the following
\begin{lemma}
  \label{winning-lines}
  In $n^d$ tic-tac-toe, there are $\displaystyle \frac{(n+2)^{d} - n^{d}}{2}$ winning lines.
\end{lemma}
Two proofs appear in~\cite{golomb-H}.  It is also useful to consider how many winning lines pass through
a given cell, as this gives an indicator of how ``powerful" the cell is.  In particular (see also
page 54 of~\cite{beck}), we have the following
\begin{lemma}
  \label{max-degree}
  Consider $n^d$ tic-tac-toe.  If $n$ is odd, then there are at most $(3^d - 1)/2$ winning lines containing
  a given cell, and this is attained only at the cell for which each coordinate is $(n+1)/2$; that is, at
  the ``center" of the board.  If $n$ is even, then there are at most $(2^d - 1)$ winning lines containing
  a given cell, with equality if there is a common $c \in [n]$ such that each coordinate $x_j$ equals either
  $c$ or $n + 1 - c$ for each $j = 1, 2, \dots, d$.
\end{lemma}

\subsection{Game play under two different types of rules}  
  
\vspace{0.5 true cm}
\noindent
We expect that most readers have played tic-tac-toe using the ``standard rules" described in the next paragraph.
For our purposes, it is helpful to recall an alternative set of rules, what Beck~\cite{beck} refers to
as ``Maker-Breaker rules".
  
\vspace{0.5 true cm}
\noindent
{\em Standard rules for the $n^d$ game:} The two players take turns claiming positions on the board.  When the
first player claims a position to occupy, they mark that position with an X; the second player marks their
occupied positions using the letter O.  Each player has the same winning objective, namely to occupy all $n$
positions along some winning line.  The first player to do so (if any) is the winner.  Here is an example of
game play in the $3^2$ game, reading from left-to-right:  
  
\vspace{0.5 true cm}
\noindent
  \bthree
      &   &   \\ \hline
      & X &   \\ \hline
      &   &
  \ethree
  \quad \quad
  \bthree
    O &   &   \\ \hline
      & X &   \\ \hline
      &   &
  \ethree
  \quad \quad
  \bthree
    O &   & X \\ \hline
      & X &   \\ \hline
      &   &
  \ethree
  \quad \quad
  \bthree
    O &   & X \\ \hline
      & X &   \\ \hline
    O &   &
  \ethree
  \quad \quad
  \bthree
    O &   & X \\ \hline
      & X & X \\ \hline
    O &   &
  \ethree
  \quad \quad
  \bthree
    O &   & X \\ \hline
    O & X & X \\ \hline
    O &   &
  \ethree
   
\vspace{0.5 true cm}
\noindent
The first player's opening move is the center square.  They make a mistake on their third move (failing to block
in the first column), allowing the second player to win the game.  This should not have happened---under the
standard rules of $3 \times 3$ tic-tac-toe, either player can force a {\em draw}, i.e., a situation in which all
nine positions in the board have been occupied without either player managing to occupy all three positions along
some winning line. Indeed, although there are three possible outcomes of this game (first player wins, second
player wins, or the game is a draw), the game will always end in a draw unless one of the players makes a
mistake.

\vspace{0.5 true cm}
\noindent
{\em Maker-Breaker rules for the $n^d$ game:} The above example game is helpful for introducing the {\em
Maker-Breaker rules} variant.  Under Maker-Breaker rules, the first player (Maker) has the same winning objective
as under standard rules:  occupying all $n$ positions along some winning line.  However, the second player
(Breaker) has a completely different winning objective:  preventing Maker from succeeding.  Whether Breaker
happens to occupy all $n$ positions along some winning line is irrelevant!  Maker wins by occupying an entire
winning line.  Breaker wins by defending every winning line against total occupation by Maker. Under Maker-Breaker
rules, in the example $3^2$ game above, Breaker has {\em not} won after their third move as they would have under
the standard rules.  In fact, Breaker is in serious trouble because Maker could win by claiming the bottom-right
position as their fourth move.

\subsection{Strategies and observations}  \label{strategies} 
  
\vspace{0.5 true cm}
\noindent
Following Appendix~C of~\cite{beck}, a {\em strategy} can be thought of as a ``handbook" that prescribes
unambiguously what a player's next move will be in every possible circumstance---a strategy is a complete plan
that tells a player exactly what to do next under any configuration of the game board.  While it is convenient for
a player to have a strategy that automates their decision-making process, strategies need not be good or sensible.
Here, we are primarily interested in whether a player has a {\em winning strategy}:  a strategy that is always
guaranteed to lead the player to achieve their winning condition before their opponent can do so.  Much of our
focus will be on {\em pairing strategies} in which the second player pairs each of their moves with the
immediately preceding move of the first player.  We now record three important observations.
    
\begin{itemize}
      
\item Observation 1:  Under Maker-Breaker rules, clearly a game cannot end in a draw.
  
\item Observation 2:  If Breaker\footnote{A word of caution:  If {\em Maker} has a winning strategy under
Maker-Breaker rules, this does {\em not} imply that the first player has a winning strategy under standard rules.
For example, Maker has a winning strategy in the $3 \times 3$ game under Maker-Breaker rules, but the first player
does not have a winning strategy under standard rules.} has a winning strategy under Maker-Breaker rules then,
under the standard rules, the second player can implement that same strategy to prevent the first player from
winning.
  
\item Observation 3:  If under standard rules the second player has a strategy (call it STR) that is guaranteed to
prevent the first player from winning, then the first player can adopt a minor adaptation of STR to prevent the
second player from winning, meaning that both players can force a draw.  Specifically, the first player can make
an arbitrary opening move which they then treat as ``invisible"---as if it never occurred.  The first player then
pretends to be the second player, and implements STR.  If STR ever requires the (true) first player to select a
position that they already occupied but regarded as invisible, they merely (i) acknowledge ownership of that
position, regarding it as visible from now on, and (ii) claim any unoccupied position on the game board and regard
it as the new invisible position.
    
\end{itemize}

\subsection{Known Results Regarding Strategies for the $n^d$ Game} \label{known-results}

\vspace{0.5 true cm}
\noindent
{\bf Two dimensions:} 
For planar ($d = 2$) tic-tac-toe under Maker-Breaker rules, Maker wins trivially if $n = 1$ or $n = 2$.  It is not
difficult to check that Maker also has a winning strategy if $n = 3$:  rotating the board if necessary, there is
no loss of generality in assuming that Maker's first two moves are the central cell $(2,2)$ and the corner  
$(3,3)$, after which they may always guarantee victory on their fourth move. (Remember that under Maker-Breaker 
rules, Maker is not at all concerned about whether Breaker happens to occupy all cells along some winning line.) 
If $n \geq 5$, it is straightforward to prove that Breaker wins by a {\em pairing strategy} (see
Figure~\ref{567}).  As
an illustration, refer to the middle panel of the figure and consider the $6^2$ game.  Breaker has a simple
winning strategy by pairing: whenever Maker chooses one of the numbered cells, Breaker immediately claims the
other cell containing the pair of that number.  Whenever Maker chooses a cell that is not numbered, Breaker
chooses any other cell that is not numbered.  Because each of the 14 winning lines contains a pair of numbers,
this strategy guarantees that Maker will occupy at least one cell on every winning line.  Breaker's winning
strategy for the $5^2$ game (left panel of Figure~\ref{567}) requires a minor modification because only one 
cell
is not numbered. Until Maker selects that cell, Breaker simply pairs moves as before:  when Maker chooses a
numbered cell, Breaker immediately chooses the other cell of the same number.  If Maker chooses the unnumbered
cell, Breaker makes an arbitrary move.  If at a later stage of the game Maker chooses a numbered cell whose pair
Breaker already occupies, Breaker simply makes another arbitrary move.  The pairing strategies in Figure~\ref{567}
can be extended to generate pairing strategies for all larger $n \times n$ boards.  For example, the pairing
strategy for the $7^2$ game in the right panel is readily obtained by encasing the $5^2$ pairing strategy with an
extra layer, and including four new pairs to cover the outermost rows and columns.

The preceding paragraph resolves the outcomes of all Maker-Breaker tic-tac-toe games on an $n^2$ board except for
$n = 4$.  ``Pre-game" pairing strategies of the sort depicted in Figure~\ref{567} cannot exist if $n = 4$, because
the 16 cells in the $4^2$ board cannot accommodate 10 distinct pairs, one pair per winning line.  Nevertheless, as
explained in Section~\ref{section44} below, Breaker still has a winning strategy by pairing,
starting with their {\em second} move.

Because Breaker has a winning strategy for the $n^2$ game if $n \geq 4$, Observations 2 and 3 above imply that,
under standard rules, either player can force a draw if $n \geq 4$.  If $n = 1$ or $n = 2$, the first player wins 
trivially under standard rules.  If $n = 3$, either player can force a draw under standard rules, and proving this
requires a small amount of case work; see pages 46--48 of~\cite{beck}.

\vspace{0.5 true cm}
\noindent
{\bf Three dimensions:} Now consider cubic ($d = 3$) tic-tac-toe, the $n^3$ game.
For $n \leq 4$, the first player has a winning strategy in the $n^3$ game under standard rules, which
naturally implies that Maker has a winning strategy under Maker-Breaker rules.  If $n = 1$ or $n = 2$, the
first player wins the $n^{3}$ game trivially.  When $n = 3$, the first player can win in four moves, the
first\footnote{Interestingly, if the $3^3$ game is played under standard rules and the first player does not open
the game by selecting the central cell $(2,2,2)$, then the second player has a winning strategy~\cite{golomb-H}.}
of which is to claim the central cell $(2,2,2)$.  By symmetry, Breaker effectively has only three possible opening
moves:  the corner $(1,1,1)$, the middle of a face $(1,2,2)$, or the cell $(1,2,1)$ in the middle of one of the   
edges of the cube.  In all three cases, Maker can select $(2,1,2)$ as their second move.  This forces Breaker to
block by selecting $(2,3,2)$ as their second move, and importantly {\em without} allowing Breaker two cells along
the same winning line.  Maker's third move is $(2,1,3)$, setting up their guaranteed win on the fourth move.
As mentioned at the beginning of this article, the $n = 4$ case was resolved by Patashnik~\cite{patashnik}:
the first player has a [very complicated] winning strategy in the $4^3$ game.

For larger $n$, Breaker can win the $n^3$ game by a pairing strategy similar to those depicted in
Figure~\ref{567}, but how large must $n$ be?  If each winning line must contain a pair, then the number of cells 
in the board (that is, $n^d$) must be at least twice as large as the number of winning lines, which is given by
Lemma~\ref{winning-lines}.  This observation, along with some routine algebra, yields the following necessary
condition for existence of a pairing strategy:
\begin{equation}
  \label{necessary}
  n \geq \frac{2}{2^{1/d} - 1}.
\end{equation}
With $d = 3$, we find that $n \geq 8$ is necessary.  In particular, when $n = 8$, there are 244 winning lines and
$8^3 = 512$ cells in the board, and $2 \times 244 \leq 512$.  See~\cite{golomb-H} for an explicit example of a
pairing strategy by which Breaker can win the $8^3$ game under Maker-Breaker rules. Generating an explicit pairing
strategy for the $9^3$ game is routine (see Section~\ref{HK}).  The explicit
pairing strategies for the $8^3$ and $9^3$ games may be used to generate explicit pairing strategies for all $n^3$
boards with $n \geq 10$, in the same way that the $5^2$ and $6^2$ games (Figure~\ref{567}) generate explicit
pairing strategies for all $n^2$ boards with $n \geq 7$.

Regarding the three remaining cases ($n = 5, 6, 7$), it has long been believed that under Maker-Breaker
rules, Breaker has a winning strategy for the $n^3$ game.  We resolve this open question in the affirmative
for the $7^3$ game.  In doing so, it is helpful to frame the problem using the language of graphs and hypergraphs.

\begin{figure}

\begin{center}
{\small
\bfive
  2 & 6  & 8 & 5 & 8  \\ \hline
 11 & 10 & 3 & 1 & 10 \\ \hline
 11 & 4  &   & 4 & 12 \\ \hline
  9 & 1  & 9 & 5 & 12 \\ \hline
  7 & 6  & 3 & 7 & 2
\efive
\qquad
\qquad
\bsix
13 &  1  &  9  &  10  &  1  &  14 \\ \hline
7  &     &  2  &  2   &     &  12 \\ \hline
3  &  8  &     &      & 11  &  3  \\ \hline
4  &  8  &     &      & 11  &  4  \\ \hline
7  &     &  5  &  5   &     &  12 \\ \hline
14 &  6  &  9  &  10  &  6  &  13
\esix
\qquad
\qquad
\bseven
    & 13 &    &    &    & 13 &    \\ \hline
15  &  2 & 6  & 8  & 5  & 8  & 16 \\ \hline
    & 11 & 10 & 3  & 1  & 10 &    \\ \hline
    & 11 & 4  &    & 4  & 12 &    \\ \hline
    &  9 & 1  & 9  & 5  & 12 &    \\ \hline
15  &  7 & 6  & 3  & 7  & 2  & 16 \\ \hline
    & 14 &    &    &    & 14 &
\eseven
}
\end{center}

\caption{Pairing strategies for the $5^2$, $6^2$ and $7^2$ games.  See text for details.}
\label{567}
\end{figure}


\section{The $4^2$ Game and Matchings} \label{section44}  

\vspace{0.5 true cm}
\noindent
The $4^2$ game is actually the source of inspiration for the proof of our main result concerning the $7^3$ game.
For $d = 2$, inequality~(\ref{necessary}) implies that $n \geq 5$ is necessary for an empty $n^2$ board to   
accommodate a pairing strategy and, in fact, we have already explained that pairing strategies really {\em do}
exist for such boards.  Unfortunately, for $n = 4$ and $d = 2$, the inequality is not satisfied.  However, in the
Maker-Breaker game on a $4^2$ board, there is still a way for Breaker to implement a pairing strategy starting
with their {\em second} move.  Following~\cite{galvin}, by rotational and reflectional symmetry, there are
effectively three opening moves that Maker can choose:  the corner $(1,1)$, the edge cell $(1,2)$, or the interior
cell $(2,2)$.  In all cases, Breaker can use $(3,3)$ as their opening move, blocking three of the ten winning
lines.  Thus, after one round of play, there are seven ``surviving" winning lines that Maker can attempt to   
occupy, and 14 unoccupied cells on the board. Because the number of unoccupied cells is twice as large as the
number of surviving winning lines and each surviving winning line contains at least two unoccupied cells, it is
feasible that Breaker could automate all of their subsequent moves by a pairing strategy. Indeed,
Figure~\ref{4-by-4} illustrates how to perform the pairing in all three cases. Maker's first move is denoted by X
and Breaker's first move is denoted by O.  From that point onward, whenever Maker chooses a numbered cell, Breaker
immediately chooses the other cell containing that number. Because each surviving winning line contains a pair,
Breaker is guaranteed to occupy at least one position on every winning line, assuring their victory.   

\begin{figure}

\begin{center}
\bfour
  7 & 5  & 1 & 1  \\ \hline
  4 & 5  & O & 6  \\ \hline
  4 & 2  & 2 & 6  \\ \hline
  X & 3  & 3 & 7
\efour
\qquad
\qquad
\bfour
  7 & 5  & 1 & 1  \\ \hline
  4 & 5  & O & 6  \\ \hline
  X & 2  & 2 & 6  \\ \hline
  4 & 3  & 3 & 7
\efour
\qquad
\qquad
\bfour
  7 & 5  & 1 & 1  \\ \hline
  4 & 5  & O & 6  \\ \hline
  2 & X  & 2 & 6  \\ \hline
  4 & 3  & 3 & 7
\efour

\end{center}

\caption{Pairing strategies for the $4^2$ games.  See text for details.}
\label{4-by-4}
\end{figure}


\subsection{Representing the $4^2$ game using graphs and hypergraphs:}

We confess that we generated the panels of Figure~\ref{4-by-4} using a bit of guess-and-check,
something that would {\bf not} be tractable when we proceed to the $7^3$ game below!  But, there
is good news:  pairing strategies like the ones in the figure can be generated {\em systematically}.  In
generating Figure~\ref{4-by-4}, for each winning line that Breaker had not yet occupied, we needed
to {\em match} a pair of cells.  The matching process can be accomplished via a straightforward
algorithm described below, and using the algorithm requires us to describe tic-tac-toe board  
configurations using the language of graphs and hypergraphs.

\vspace{0.5 true cm}
\noindent
A finite {\em hypergraph} $\HH = (V(\HH), E(\HH))$ is a pair consisting of two sets:  $V(\HH)$ is a finite set
whose elements are called {\em vertices}; $E(\HH)$ is a subset of the power set of $V(\HH)$, and elements of
$E(\HH)$ are called {\em edges}.  A hypergraph is called {\em $k$-uniform} if all edges have cardinality $k$.   
(A $2$-uniform hypergraph is a simple, undirected {\em graph}.)
The {\em degree} of $v \in V(\HH)$ is the number of edges in $E(\HH)$ that contain $v$. A hypergraph is called
{\em almost disjoint} (or {\em linear} in some references) if given any distinct edges $e_{1}, e_{2} \in E(\HH)$,
we have $|e_{1} \cap e_{2}| \leq 1$.  (That is, a hypergraph is almost disjoint if any two distinct edges
intersect in at most one vertex.)

The empty $4^2$ tic-tac-toe board is associated with a hypergraph $\Hempty$ with 16 vertices and 10 edges, where
vertices correspond to cells in the board and edges correspond to winning lines.  This hypergraph is $4$-uniform
because every winning line consists of precisely 4 cells, and is almost disjoint because two distinct winning
lines cannot have more than one cell in common.  By Lemma~\ref{max-degree}, the maximum degree of any vertex of
$\Hempty$ is 3.

The pairing strategies in Figure~\ref{4-by-4} were generated by considering possible configurations of the
board after Maker and Breaker had each made one move, such that Breaker's first move eliminated 3 winning lines.
Here is a worked example of how to associate a {\em non-uniform}, almost disjoint hypergraph $\HH$ that
represents the state of the $4^2$ game board after each player has made one move.  Suppose that the game
board looks like this
\begin{center}
\bfour
  &   &     & \\ \hline
  & X &     & \\ \hline
  &   &  O  & \\ \hline
  &   &     &
\efour
\end{center}
where $X$ denotes Maker's first move and $O$ denotes Breaker's first move.  Now, number the 14 remaining
unoccupied cells:

\begin{center}
\bfour
1  & 2  &  3  & 4  \\ \hline
5  & X  &  6  & 7  \\ \hline
8  & 9  &  O  & 10 \\ \hline
11 & 12 & 13  & 14
\efour
\end{center}
The 7 surviving winning lines give rise to a hypergraph $\HH$ with vertex set $V(\HH) = [14]$ and
whose edge set $E(\HH)$ contains 7 edges
$$
\{1,2,3,4\} \; \; \; \{5,6,7\} \; \; \; \{4,6,9,11\} \; \; \; \{2,9,12\}
$$
$$
\{1,5,8,11\} \; \; \;  \{4,7,10,14\} \; \; \;  \{11,12,13,14\}
$$
which we denote as $e_{1}$ through $e_{7}$.

To set up the matching problem, consider the bipartite graph $G = (X(G), Y(G), E(G))$ whose
vertex set $V(G) = X(G) \sqcup Y(G)$ consists of partite classes
$$
  X(G) = \{ e_{1}^{1}, e_{1}^{2}, e_{2}^{1}, e_{2}^{2}, \dots, e_{7}^{1}, e_{7}^{2} \}, \qquad
  Y(G) = \{1, 2, 3, \dots, 14\}
$$
and for each $j = 1, 2, \dots, 7$ we have $e_{j}^{1} = e_{j}^{2} = e_{j}$.  A vertex
in $e_{j}^{i} \in X(G)$ is connected to $m \in Y(G)$ if and only if $m$ is an element
of the set $e_{j}^{i}$.  Note that elements of the partite class $X(G)$ are
edges from $E(\HH)$, with each edge listed twice.

\begin{figure}
 \begin{center}
   \includegraphics[width=0.3\textwidth]{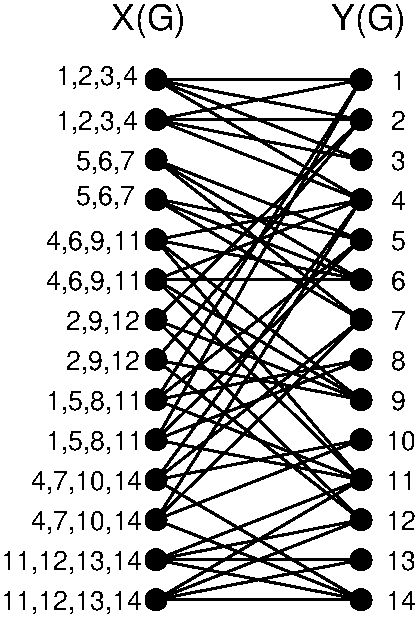}
   \qquad
   \includegraphics[width=0.3\textwidth]{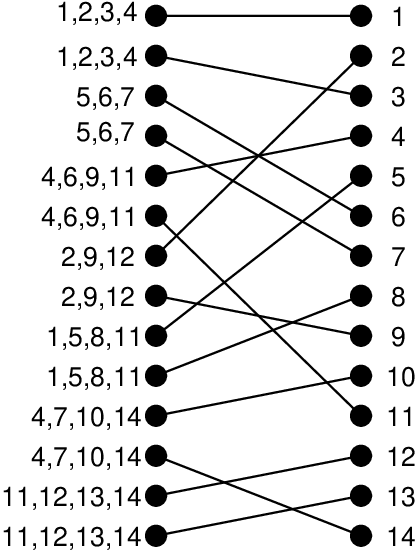}
 \end{center}
 \caption{Left panel:  The bipartite graph $G = (X(G), Y(G), E(G))$ associated with the
 $4^2$ board described in the worked example in the main text.  Right panel:  An example of a maximal
 matching of size 14 from the graph in the left panel.}
 \label{HK44}
\end{figure}

The graph $G$ is shown in the left panel of Figure~\ref{HK44}.  Remember that each vertex in $X(G)$ corresponds to
a list of unoccupied vertices in some winning line that Breaker has not yet occupied.  Each vertex in $Y(G)$
represents an unoccupied cell in the game board.  The right-panel of Figure~\ref{HK44} shows an example of an {\em
$X(G)$-perfect matching}:  a set of pairwise non-adjacent edges of $G$ such that each vertex in $X(G)$ is
contained in one of the edges. Associated with an $X(G)$-perfect matching is an injective function $f: X(G) \to
Y(G)$ such that for each $e_{j}^{(i)} \in X(G)$, the pair $(e_{j}^{(i)}, f(e_{j}^{(i)}))$ is in $E(G)$. Having an
$X(G)$-perfect matching generates a pairing strategy by which Breaker can win. For each $j \in \{1, 2, \dots,
7\}$, note that $f(e_{j}^{(1)})$ and $f(e_{j}^{(2)})$ are distinct elements of $Y(G)$, and they correspond to   
distinct cells on some winning line that Breaker has not yet occupied.  Breaker's pairing strategy is to pair
these cells:  henceforth, whenever Maker selects an unoccupied cell $f(e_{j}^{(1)})$ or $f(e_{j}^{(2)})$, Breaker
immediately counters by selecting the other member of that pair, unless they already occupy it. For example,
referring back to our earlier numbering of the cells in the game board, the right panel of Figure~\ref{HK44}
generates the pairing strategy   
\begin{center}
\bfour
1  & 2  &  3  & 4  \\ \hline
5  & X  &  6  & 7  \\ \hline
8  & 9  &  O  & 10 \\ \hline
11 & 12 & 13  & 14
\efour
\qquad
$\longrightarrow$
\qquad
\bfour
a  & e  &  a  & g  \\ \hline
d  & X  &  b  & b  \\ \hline
d  & e  &  O  & f  \\ \hline
g  & c  &  c  & f
\efour
\end{center}
where paired cells are indicated using the letters $a$ through $g$.

{\em Remark:} Another way to think of finding a maximal matching in the bipartite graph $G$ is that we need to
produce a set of distinct representatives from the 14 sets associated with vertices in partite class $X(G)$.

\subsection{Finding Matchings:  The Hopcroft-Karp Algorithm}  \label{HK}

One last question remains:  Given a bipartite graph like the one in the left panel of Figure~\ref{HK44},
how can we {\em find} and $X(G)$-perfect matching like the one in the right panel?

\vspace{0.5 true cm}
\noindent
The Hopcroft-Karp algorithm~\cite{hopcroft-K} generates a maximal matching in a finite bipartite graph, i.e., a
maximum cardinality set of edges for which no two edges share a common vertex.  See~\cite{hopcroft-K} for a
thorough explanation of how the algorithm works.  Like some its earlier predecessors, the Hopcroft-Karp algorithm
begins with a partial matching, and attempts to find a larger partial matching by searching for an {\em augmenting
path}.  Let us recall how this works. Suppose that $G = (X(G), Y(G), E(G))$ is a bipartite graph.  The vertex set
$V(G) = X(G) \sqcup Y(G)$ is a disjoint union of nonempty partite classes $X(G)$ and $Y(G)$, and each edge in
$E(G)$ consists of precisely one vertex from $X(G)$ and one vertex from $Y(G)$.  To avoid triviality, assume
that $E(G)$ is nonempty.  Now,
\begin{itemize}
  \item Begin with any nonempty partial matching $P \subseteq E(G)$.  Here, $P$ could consist of a single edge.   
        If $P$ contains more than one edge, no two edges of $P$ share a vertex in common.
  \item Attempt to find an {\em augmenting path}:  A path beginning at an unmatched vertex (one that does not 
appear
        in any edge in $P$), ending at an unmatched vertex, and following a sequence of edges $e_1, e_2, \dots   
        e_{2k + 1}$ for which $e_j \in P$ if $j$ is even and $e_{j} \notin P$ if $j$ is odd.  Here, $k$ is a 
        non-negative integer; $k = 0$ is allowed.  If such an augmenting path exists, update $P$ by removing    
        $e_2, e_4, \dots e_{2k}$ and adding $e_{1}, e_{3}, \dots, e_{2k+1}$.  Doing so increases the cardinality
        of $P$ by 1, resulting in a larger partial matching.
  \item Repeat the preceding step for as long as there exists an augmenting path.  When the process terminates,
        $P$ is a matching of maximal size.
\end{itemize}


\section{Main result:  Resolving the $7^3$ game}   \label{section-main}

\vspace{0.5 true cm}
\noindent
In this section, we focus exclusively on the $7^3$ game under Maker-Breaker rules.  The game board has $7^3 = 343$
cells and, by Lemma~\ref{winning-lines}, there are 193 winning lines. Since $193 \times 2 > 343$, Breaker does not
have a pairing strategy prior to the start of game play---as explained above, such strategies exist for the $n^3$
game if and only if $n \geq 8$.  Instead, we will reuse the idea from the above proof that Breaker wins the $4^2$
game:  if Breaker's initial moves are sufficiently effective in blocking winning lines, then a pairing strategy
always becomes available after several rounds of game play.  Finding such pairing strategies can be accomplished
by (i) associating graphs and hypergraphs with game board configurations, precisely as in the $4^2$ example above,
and (ii) implementing the Hopcroft-Karp matching algorithm.




\subsection{Main result}   

\vspace{0.5 true cm}
\noindent
Now, we state and prove our main result:
\begin{theorem}  \label{main-theorem}
Under Maker-Breaker rules, Breaker has a winning strategy in $7^3$ tic-tac-toe.
\end{theorem}
The following corollary follows immediately from the Observations 2 and 3 above.
\begin{corollary}
  Under standard rules, either player can force a draw in $7^3$ tic-tac-toe.
\end{corollary}
The key idea of the proof is to mimic what was done previously for the $4^2$ game:
Breaker will use their first few moves to occupy as many winning lines as possible,
after which they are able to automate all subsequent moves by a simple pairing strategy.  We will consider two
cases according to whether Maker's first move is the central cell $(4,4,4)$.  In what follows, at any stage of 
game
play, we say that a winning line is a {\em survivor} if Breaker does not yet occupy any cells on that line;
otherwise we say that Breaker has {\em eliminated} the winning line.  We also single out four of the winning lines
that we refer to as {\em superdiagonals}:  winning lines that run between opposite corners of the cubic game      
board, and are not parallel to any face of the cube.  More precisely, the coordinates of the cells on the four
superdiagonals are
        
\begin{center}
\begin{tabular}{llllll}
  superdiagonal {\bf (A)}:   & (1,1,1) & (2,2,2) & $\cdots$ & (6,6,6)  & (7,7,7)   \\
  superdiagonal {\bf (B)}:   & (1,7,1) & (2,6,2) & $\cdots$ & (6,2,6)  & (7,1,7)   \\
  superdiagonal {\bf (C)}:   & (7,1,1) & (6,2,2) & $\cdots$ & (2,6,6)  & (1,7,7)   \\
  superdiagonal {\bf (D)}:   & (7,7,1) & (6,6,2) & $\cdots$ & (2,2,6)  & (1,1,7).
\end{tabular}
\end{center}

\subsection{Proof of Theorem~\ref{main-theorem} if Maker's first move is $(4,4,4)$:}  \label{444}

In the hypergraph $\Hempty$ described above, the vertex
corresponding to the cell $(4,4,4)$ has degree 13, and is the only vertex with a degree higher than 7.  The
vertices of degree $7$ are the ones corresponding to the other cells on the superdiagonals.  Breaker will
use their first four moves to occupy one cell on each superdiagonal in such a way that they eliminate
27 winning lines.  By doing so, after each player has made exactly four moves, 166 of the winning lines
are survivors.  At this stage of the game, the board has 335 unoccupied cells and, since
$166 \times 2 < 335$, it is feasible that there may exist a pairing strategy that Breaker can exploit for
the remainder of the game.  Indeed, we will show that this always happens.

\vspace{0.5 true cm}
\noindent
{\em Step 1:  Selecting Breaker's first four moves.}
Maker's first move is $(4,4,4)$.  Breaker chooses $(5,5,5)$
on superdiagonal {\bf (A)} as their first move, thereby eliminating 7 winning lines.  In preparation for
computer-assisted case work, there is one symmetry\footnote{We do this mainly as an acknowledgment
that there are ways to leverage symmetry to reduce future casework.  Admittedly, a more careful consideration of
symmetries of the cube could have simplified our computational efforts, but the computer-assisted portion of our
work was fast enough that we did not feel the need.} we now exploit.  Both players' first moves lie on 
superdiagonal
{\bf (A)}.  By rotating the cube about superdiagonal {\bf (A)} if
necessary, we may without loss of generality assume that Maker's second move is of the form $(x,y,z)$ where $x
\leq y$.  This ensures that the cell $(6,2,6)$ on superdiagonal {\bf (B)} is still available.  Breaker's second
move is to occupy $(6,2,6)$, which eliminates 7 winning lines that had not been eliminated previously.
Regardless of Maker's third move, at least two of the corner\footnote{That is, cells for which each coordinate   
is either $1$ or $7$.} cells on superdiagonals {\bf (C)} or {\bf (D)} remain unoccupied.
Breaker's third move is to choose such a corner cell lying on one of these two superdiagonals,
thereby eliminating 7 winning lines that had not been eliminated previously.  Finally, regardless of Maker's 
fourth move, Breaker's fourth move is to choose any unoccupied cell along whichever superdiagonal, either
{\bf (C)} or {\bf (D)}, they have not yet occupied.  In doing so, Breaker's fourth move eliminates 6 winning
lines that had not been eliminated previously---not 7, because this move will lie along one winning line that
Breaker had previously eliminated.  In summary, after each player has completed four moves, Breaker has
eliminated 27 winning lines and there are 166 survivors.

{\em Remarks:}  (i)  Breaker's third and fourth moves can be prescribed in such a way that their choices are
fully determined.  For instance, Breaker's third move is to select the first unoccupied cell among the
ordered list $(7,1,1)$, $(1,7,7)$, $(7,7,1)$, and $(1,1,7)$.  In doing so, the number of possible configurations
of the game board after Breaker's fourth move is completely determined by the number of options that Maker
had.  Maker had 194 options for the second move, the number of cells of the form $(x,y,z)$ with $x \leq y$
that remain unoccupied after Breaker's first move.  Maker's third move is chosen from among 339 unoccupied cells,
and their fourth move from among 337.  This means that the number of possible board configurations at this
stage is
$$
  194 \times 339 \times 337 = 22,163,142.
$$
(ii) Because each player has made only 4 moves, each surviving winning line contains at least 3 unoccupied cells
and Breaker is not in imminent danger of losing the game.

\vspace{0.5 true cm}
\noindent
{\em Step 2:  Framing the search for a pairing strategy as a matching problem.}  Given a configuration of
the game board after Breaker's fourth move and presuming that the procedure in Step 1 has been followed,
we associate a hypergraph $\HH = (V(\HH),E(\HH))$ as follows.  The vertex set $V(\HH)$ consists of the 335 cells
that are not yet occupied.  The edge set $E(\HH)$ consists of 166 edges that correspond to the survivor winning
lines.  Each edge consists of all {\em unoccupied} vertices of a survivor.  The hypergraph $\HH$
is {\bf not} 7-uniform because there are survivors for which Maker occupies at least one cell.  Every edge
in $E(\HH)$ contains between 3 and 7 vertices.  (In order for there to exist an edge that contains precisely 3
vertices, all four of Maker's opening moves would have to lie along a single winning line, in which case the
edge corresponding to that winning line is the only one containing only 3 vertices.) 
  
To set up the matching problem, let $e_{1}, e_{2}, \dots, e_{166}$ be the edges in $E(\HH)$.
Consider the bipartite graph $G = (X(G), Y(G), E(G))$ whose vertex set $V(G) = X(G) \sqcup Y(G)$ consists of 
partite
classes
$$
  X(G) = \{ e_{1}^{1}, e_{1}^{2}, e_{2}^{1}, e_{2}^{2}, \dots, e_{166}^{1}, e_{166}^{2} \}, \qquad
  Y(G) = \{1, 2, 3, \dots, 335\}
$$
and for each $j = 1, 2, \dots, 166$ we have $e_{j}^{1} = e_{j}^{2} = e_{j}$.  A vertex
in $e_{j}^{i} \in X(G)$ is connected to $m \in Y(G)$ if and only if $m$ is an element
of the set $e_{j}^{i}$.  Note that {\em vertices} in the partite class $X(G)$ are
{\em edges} from $E(\HH)$, with each edge duplicated twice, and vertices in the partite class $Y(G)$ 
correspond to the unoccupied cells in the tic-tac-toe board.  If we can find an $X(G)$-perfect matching,
then we can generate a pairing strategy that Breaker can use to win the game.  More precisely, suppose that
there exists an injective function $f: X(G) \to Y(G)$.  For each $j \in [166]$, the cells
$f(e_{j}^{(1)})$ and $f(e_{j}^{(2)})$ are distinct, and we pair them with each other.  In doing so, we
create 166 pairs among the unoccupied cells in such a way that each surviving winning line contains one
pair.  Three of the unoccupied cells in the game board are not paired.  Breaker's strategy for the remainder
of the game is exactly like the strategy for the $5^2$ game.  Whenever Maker selects an unoccupied cell
$f(e_{j}^{(1)})$ or $f(e_{j}^{(2)})$, Breaker immediately counters by choosing the other member of that pair
unless they already occupy it.  If Maker ever selects one of the three unpaired cells, or a cell whose pair
Breaker already occupies, Breaker's next move is arbitrary.  This strategy guarantees that Breaker will 
occupy at least one cell on every winning line, assuring their victory.

It remains to show that the pairing described in the preceding paragraph always exists, for every possible
configuration of the board after Step 1 was completed.

\vspace{0.5 true cm}
\noindent
{\em Step 3:  Finding matchings to generate pairing strategies for Breaker.}  The final step of the proof was    
computer-assisted, and our code and compilation instructions is available here:
\begin{center}
  {\color{blue} \url{https://people.math.harvard.edu/~jcain2/tic-tac-toe-code.tar.gz}}
\end{center}
For each of the 22,163,142 possible board configurations described in Step 1, we constructed the associated  
bipartite graph $G$ described in Step 2.  Then, we used the Hopcroft-Karp matching algorithm (see
reference~\cite{hopcroft-K} or below) to determine the largest size of a matching between the
partite classes $X(G)$ and $Y(G)$.  In every case, our program determined that the maximum matching size is $332 =
2 \times 166$, indicating that a pairing strategy {\em always} fits.  Our code actually constructs the maximal
matchings, thereby creating the pairing strategies.

\subsection{Proof of Theorem~\ref{main-theorem} if Maker's first move is not $(4,4,4)$:}  \label{not444}

The previous subsection considered only one of {\em many} possible opening moves that Maker might select, and the
reader might [rightfully] worry that far more casework remains.  Luckily, if Maker's first move is not $(4,4,4)$,
less casework is involved because Breaker can establish a viable pairing strategy after their {\em third} move.  
Here is a modification of Step 1 from the preceding subsection; Steps 2 and 3 proceed in the same manner as 
before.

\vspace{0.5 true cm}
\noindent
{\em Step 1:  Selecting Breaker's first three moves.}
Maker's first move is not $(4,4,4)$, so Breaker takes $(4,4,4)$ as their first move, eliminating
13 winning lines.  Rotating the cube if necessary, there is no loss of generality in assuming that
Maker's first {\em and} second moves are in the ``half cube" for which the first coordinate is at
most $4$.  Breaker selects $(5,5,5)$ on superdiagonal {\bf (A)} as their second move, eliminating 6
winning lines that were not eliminated by their first move.  Regardless of Maker's third move, at
least one of the cells $(6,2,6)$ or $(7,1,7)$ on superdiagonal {\bf (B)} is unoccupied because
Maker's first two moves cannot have $6$ or $7$ in their first coordinate.  Breaker's third move
is to choose one of these two cells; if both are unoccupied, then for the sake of definiteness
Breaker chooses $(6,2,6)$.  In any case, Breaker's third move eliminates 6 winning lines that they
had not eliminated previously.  In summary, after Breaker's third move, 25 of the original winning
lines have been eliminated, leaving 168 survivors.  There are now 337 unoccupied cells remaining on
the board and, since $168 \times 2 = 336 < 337$, there is hope that the unoccupied cells can
accommodate 168 pairs, with one pair per surviving winning line.

\vspace{0.5 true cm}
\noindent
{\em Steps 2 and 3} are performed exactly as in the preceding subsection.  Because Breaker's moves are
completely prescribed, the number of cases to check is based upon how many choices Maker has for
each of their first three moves.  Maker's first move $(x_1, y_1, z_1)$ has 195 possibilities since
$x_1 \leq 4$ and they do not choose $(4,4,4)$.  Maker's second move $(x_2, y_2, z_2)$ has 194 possibilities
because $x_2 \leq 4$
and both Maker and Breaker selected cells with first coordinate at most 4 during the first round of
play.  Maker's third move has 339 options, yielding
$$
  195 \times 194 \times 339 = 12,824,370
$$
possible board configurations after three rounds of play.  In every case, the Hopcroft-Karp algorithm
succeeded in producing a maximal matching of size $336 = 2 \times 168$ in the bipartite graph
generated from the game board.  Thus, after their third move, Breaker is able to win the game via
a simple pairing strategy.  The proof of Theorem~\ref{main-theorem} is now complete.  $\square$

\vspace{0.5 true cm}
For the record, the total time required for the first author's laptop computer to perform all
casework from Step 3 of the proof was approximately 5 hours.
Had we exercised more effort to reduce redundancy, the computational time could have been reduced to
less than an hour.  Doing so would have resulted in a less concise explanation of Step 1 in the proof,
something we preferred to avoid.


\section{A non-constructive approach}  \label{section-bigamy}

The proof of Theorem~\ref{main-theorem} involved construction of matchings via application of the Hopcroft-Karp 
algorithm.  Now, we present a different, non-constructive proof that Breaker can win the $4^2$ Maker-Breaker game.  
Once again, we will frame the game using the language of hypergraphs.  By invoking a variant of the matching 
theorem of K\H{o}nig, Hall, and Egerv\'{a}ry, we will show that Breaker can arrange to have a pairing strategy 
starting from their second move.  We will also report on our attempt to use this approach to prove that Breaker 
can win the $7^3$ game.

Suppose that Maker and Breaker have each made one move in the $4^2$ game, with Breaker's move eliminating three 
winning lines. As in Section~\ref{HK}, let $\HH = (V(\HH), E(\HH))$ be the hypergraph associated with the game 
board:  the 14 vertices in $V(\HH)$ represent the unoccupied cells, and each of the 7 edges in $E(\HH)$ contains 
the unoccupied positions of a surviving winning line. The hypergraph is not $4$-uniform, but we know that each 
edge contains either three or four vertices.  The following corollary of the K\H{o}nig-Hall-Egerv\'{a}ry matching 
theorem appears as Theorem 1.4 in~\cite{beck2002}.  In reference~\cite{beck}, Beck refers to this result as the 
``bigamy corollary" of Hall's theorem.

\begin{proposition}  \label{bigamy-prop}
 If for every $\mathcal{G} \subseteq E(\HH)$ the inequality
 \begin{equation}
   \label{bigamy-ineq}
   \left| \bigcup_{S \in \mathcal{G}} S \right| \geq 2 |\mathcal{G}|
 \end{equation}
 holds, then Breaker can win the game by a pairing strategy.
\end{proposition}

Proposition~\ref{bigamy-prop} can be used to give a non-constructive proof of a result we have already stated:
\begin{proposition}  \label{44again}
  In $4^2$ Maker-Breaker tic-tac-toe, if Breaker's first move eliminates 3 winning lines, then Breaker has a winning
  strategy by pairing for all of their remaining moves.
\end{proposition}

{\em Sketch of proof:}  Without loss of generality, we may assume that the first moves of Maker and Breaker fall into one of
the three cases in Figure~\ref{4-by-4}.  To apply Proposition~\ref{bigamy-prop}, we check that inequality~(\ref{bigamy-ineq})
always holds.  Let $m = |\mathcal{G}|$.  Applying the inclusion-exclusion principle to estimate the cardinality of 
the union on the left-hand side yields
\begin{equation}
  \label{IE}
   \left| \bigcup_{S \in \mathcal{G}} S \right| \; = \; \sum_{i = 1}^{m} |S_{i}| \; - \; 
   \sum_{1 \leq i < j \leq m} | S_{i} \cap S_{j}| \; + \; \sum_{1 \leq i < j < k \leq m} | S_{i} \cap S_{j} \cap S_{k} |.
\end{equation}
Notice that this expansion terminates---there is no need to consider intersections of more than three elements of $\mathcal{G}$
because in $4^2$ tic-tac-toe there are never more than three winning lines that pass through a given cell on the board.
At this stage, one may estimate each term on the right-hand side of~(\ref{IE}) to argue that regardless of $\mathcal{G}$,
inequality~(\ref{bigamy-ineq}) is satisfied.  The crude estimate
\begin{equation}
  \label{crude}
   \left| \bigcup_{S \in \mathcal{G}} S \right| \; \geq \; \sum_{i = 1}^{m} |S_{i}| \; - \; 
   \sum_{1 \leq i < j \leq m} | S_{i} \cap S_{j}| \; \geq 3m - \binom{m}{2}
\end{equation}
follows from the fact that each edge in $E(\HH)$ contains at least three vertices, and $|S_{i} \cap S_{j}| \leq 1$ 
due to the fact that $\HH$ is an almost disjoint hypergraph.  Inequality~(\ref{crude}) shows 
that~(\ref{bigamy-ineq}) holds if $m = |\mathcal{G}|$ is $1$, $2$, or $3$.  The largest $m$ that we must consider 
is $7$, but~(\ref{bigamy-ineq}) certainly holds (with equality) in that case.  With a small amount of additional 
care than was used in the crude estimate~(\ref{crude}), it is straightforward to check that~(\ref{bigamy-ineq}) 
holds if $m = 4,5,6$.  $\square$

Incidentally, when estimating the terms in the inclusion-exclusion expansion, we found it easier to (i) 
temporarily ignore Maker's first move, so that $\HH$ could be regarded as a $4$-uniform, almost disjoint 
hypergraph; (ii) estimate the terms in the inclusion-exclusion expansion to obtain an {\em incorrect} lower bound 
for $\left| \bigcup_{S \in \mathcal{G}} S \right|$; and (iii) explain how to correct the lower bound after 
reinstating Maker's first move.

We attempted to apply Proposition~\ref{bigamy-prop} to the $7^3$ Maker-Breaker game, to generate a 
non-constructive proof that Breaker can win.  As before, we presume that each player has completed four moves, 
during which Breaker has eliminated 27 winning lines, leaving 166 survivors.  The hypergraph $\HH$ associated with 
this stage of the game has 335 vertices corresponding to unoccupied cells of the board, and 166 edges, each 
containing the unoccupied cells on some surviving winning line. As one would expect, showing that 
inequality~(\ref{bigamy-ineq}) is satisfied for every choice of $\mathcal{G} \subseteq E(\HH)$ is challenging.  
Routine estimates of terms in an inclusion-exclusion expansion reveal that the inequality holds if $|\mathcal{G}|$ 
is suitably small, say 20 or less.  Likewise, if $|\mathcal{G}|$ is suitably large (close to 166), it is easy to 
check that the inequality holds.  There is an interesting duality\footnote{In a sense that can be made precise 
using the language of dual hypergraphs} between the process of verifying that~(\ref{bigamy-ineq}) holds for small 
$|\mathcal{G}|$ versus large $|\mathcal{G}|$. We will not bother to elaborate, because our attempt to apply 
Proposition~\ref{bigamy-prop} failed miserably for a wide range of intermediate $|\mathcal{G}|$.

\section{Discussion and further explorations}  \label{section-discussion}

In this article, we proved that Breaker has a winning strategy in $7^3$ Maker-Breaker tic-tac-toe, from which it 
follows that either player can force a draw in the $7^3$ game under standard rules.  We also explained how to 
frame the process of finding pairing strategies using the language of matching theory.  Before concluding, we 
would be remiss not to discuss (i) an alternative approach involving potential functions for proving that Breaker 
can win certain $n^d$ games and (ii) the only two remaining open cases for the $n^3$ Maker-Breaker game:  $n = 5$ 
and $n = 6$.

\subsection{Using potential functions to measure Breaker's danger}    \label{potential}

The \ES theorem~\cite{erdos-S} (see also Section 10 of~\cite{beck}) gives sufficient conditions for Breaker to 
have a pairing strategy that leverages a {\em potential function}.  Consider the empty board for the $n^d$ 
Maker-Breaker game, and let $\HH = (V(\HH), E(\HH))$ be the associated $n$-uniform, almost disjoint hypergraph.  
Let $\Delta$ denote the maximum degree among all $v \in V(\HH)$.
\begin{theorem}
  \label{EST}
  (\ES theorem)  If $\Delta + |E(\HH)| < 2^{n}$, then Breaker has a winning strategy.
\end{theorem}
For instance, in the $8^3$ game, we have $n = 8$, $\Delta = 7$, and the number of winning lines is $|E(\HH)| = 244$.
Since $7 + 244 < 2^{8} = 256$, Breaker has a winning strategy for the $8^3$ game.  

The proof of Theorem~\ref{EST}, which is short and readable, makes use of a potential function that measures the 
{\em total danger} that Breaker faces at each stage of the game.  Each surviving winning line presents a danger to 
Breaker.  If a surviving winning line contains $u$ unoccupied cells, it is assigned a danger of $2^{-u}$.  The 
total danger is the sum of the dangers of all surviving winning lines.  Winning lines that Breaker has eliminated 
have zero danger.  Each time Maker chooses a cell, they double the danger of every survivor containing that cell. 
For each $i = 0, 1, 2, \dots$, let $S_{i}$ denote the set of all surviving winning lines immediately after Maker's 
$i$th move and, for $s \in S_{i}$, let $u_{s}$ denote the number of unoccupied cells on $s$.  Then the total 
danger that Breaker faces immediately after Maker's $i$th move is
$$
  D_{i} = \sum_{s \in S_{i}} 2^{-u_s}.
$$
For example, in the $7^3$ game before Maker's first move, there are 193 survivors each with an associated danger 
of $2^{-7}$, for a total danger of $D_{0} = 193/128$.

\begin{figure}
  \begin{center}
   \includegraphics[width=0.45\textwidth]{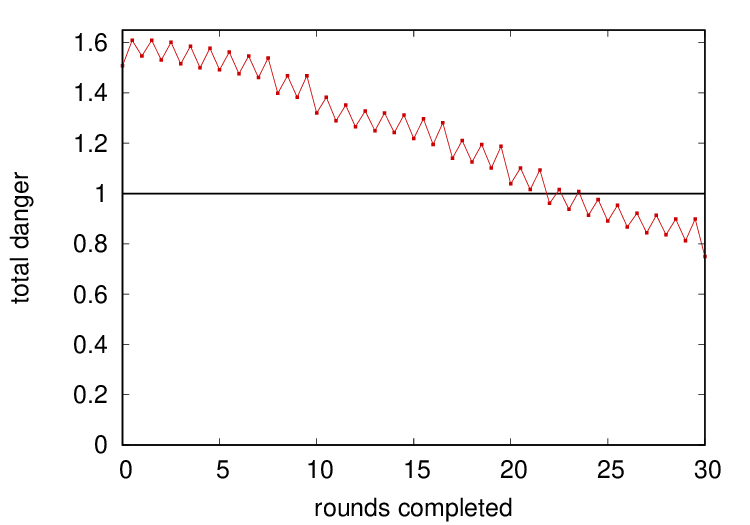}
   \includegraphics[width=0.45\textwidth]{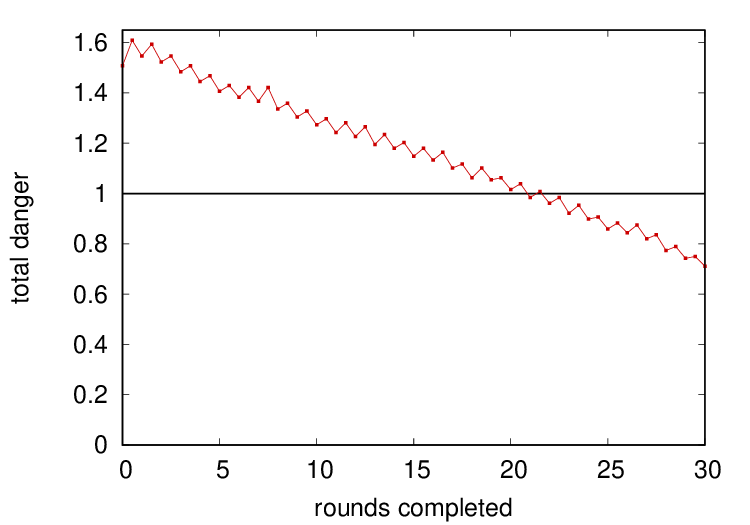}
  \end{center}
  \caption{Left:  Total danger after each move, assuming that Maker always chooses moves that contribute the
  highest possible increase to total danger and Breaker always chooses moves that contribute the highest
  possible decrease to total danger.  Right:  Breaker uses the same strategy as in the left panel, but
  Maker makes random moves (apart from their first move, which is the center cell).}
  \label{fig-danger}
\end{figure}

One consequence of the proof Theorem~\ref{EST} is that, if Breaker ever manages to arrange that $D_i < 1$ for some 
$i \in \mathbb{N}$, they have a very simple winning strategy:  always select a cell that achieves a maximal 
reduction in danger.  We considered attempting to prove that Breaker can always achieve $D_{i} < 1$ in the $7^3$ 
game, but we found that the number of rounds $i$ needed to accomplish this can be as high as 24.  For example, 
Figure~\ref{fig-danger} shows what happens if Maker adopts a greedy strategy, always selecting a cell that 
contributes the highest possible increase to total danger.  Meanwhile, Breaker always counters with moves that 
decrease total danger as much as possible.  The first $i$ for which total danger $D_{i} < 1$ is $i = 24$.  The 
right panel of the figure shows what happens if Maker's first move is the center cell, but the rest of their moves 
are chosen randomly.


\subsection{The $5^3$ and $6^3$ games}

Now, $n = 5$ and $n = 6$ are the only unresolved cases for which it is unknown whether Breaker has a winning 
strategy in Maker-Breaker tic-tac-toe on an $n^3$ board.  It is believed that in both of these cases, Breaker has 
a winning strategy.  Some evidence in support of this belief can be obtained by reproducing the left panel of 
Figure~\ref{fig-danger} for the $5^3$ and $6^3$ games.  In both the $5^3$ and $6^3$ games, if Maker follows the 
same greedy approach described in the preceding section, we find that Breaker can arrange for $D_{i} < 1$ when $i 
= 31$.

What if we attempt to replicate the matching-based approach that was used to resolve the $7^3$ game? Unfortunately, 
this is not feasible for the $5^3$ game, which has 109 winning lines but only 125 cells.  Breaker would require 
{\em many} rounds of game play to eliminate enough winning lines to establish viability of a pairing strategy.  
By the time Breaker manages to do so, the number of possible board configurations is staggering, and the 
task of testing each such configuration for a pairing strategy is computationally intractable.  The challenge is 
further exacerbated by the fact that Breaker cannot replicate the simple early-game strategy from Step 1 of the 
proof of Theorem~\ref{main-theorem}.  Merely eliminating winning lines as quickly as possible will not suffice.
Breaker's strategy must also account for ``emergency" scenarios, i.e., they must know how to react if Maker 
establishes a configuration that presents a serious threat.  Chapter VII of~\cite{beck} offers a detailed 
discussion of how Breaker might deal with emergencies.

The $6^3$ game is not quite as hopeless.  The game board has 148 winning lines and 216 cells.  Unlike the $n^3$ 
boards with $n \geq 3$ odd, there is no central cell of degree 13.  The largest number of winning lines passing 
through a cell is 7, and the four superdiagonals do not intersect one another.  As in the $7^3$ game, in their 
first four moves, Breaker can arrange to occupy at least one cell on each superdiagonal in such a way that 27 
winning lines are eliminated.  After both players have made four moves, this means that there are 121 survivors 
and 208 unoccupied cells, and Breaker is not in imminent danger of losing on the next move.  It is possible that 
Breaker could, within the next few rounds of play, eliminate enough additional winning lines to establish a viable 
pairing strategy, while making sure to block any potential threats from Maker.  With sufficient care beyond the 
scope of what we are willing to attempt, it may be possible to resolve the $6^3$ game computationally.

Finally, let us express our hopes and encouragement that some reader(s) may use the ideas in this article to prove 
Theorem~\ref{main-theorem} without computer assistance to handle casework.  Because we authors have very limited 
experience in combinatorics, we may be woefully unaware of theorems that could have elicited elegant proofs of our 
main result.

\newpage
\begin{center}
  {\bf Acknowledgments}
\end{center}

JWC gratefully acknowledges the research support of the Budapest Semesters in Mathematics Director's Mathematician 
in Residence Program.  He also expresses gratitude to the actor Matthew Broderick, whose role in the movie {\em 
War Games} served as a source of inspiration to write a program automating a computer to play tic-tac-toe against 
itself to generate Figure~\ref{fig-danger}.  IMR and NCK gratefully acknowledge financial support from the Harvard 
College Research Program and the Center of Mathematical Sciences and Applications.


\end{document}